\documentclass[12pt]{article}
\usepackage{latexsym}
\usepackage{amsmath}
\usepackage{amsthm}
\usepackage{enumerate}
\voffset -.55in\newtheorem{theorem}{Theorem}[section] 
\newtheorem{lemma}{Lemma}[section]
\newtheorem{definition}{Definition}[section]

\newtheorem{remark}{Remark}[section]
\newtheorem{proposition}{Proposition}[section]

\newcommand{\R}{{\bf R}}

\newcommand{\Id}{\mbox{\it Id}}

\title{Two dimensional compact simple Riemannian manifolds are
boundary distance rigid}
\author{Leonid Pestov\thanks{Part of this work was done while the author\/ was visiting MSRI and
the University of Washington.}\\
Institute of Computational Mathematics and Mathematical Geophysics\\Russian
Academy of Sciences\\Novisibirsk 630090, Russia
 \and Gunther Uhlmann\thanks{Partly supported by NSF and a John
Simon Guggenheim Fellowship.}\\
Department of Mathematics\\University of
Washington\\Seattle, WA 98195, USA}
\date{}
\begin{document}
\renewcommand{\theequation}{\thesection.\arabic{equation}}
\maketitle

\begin{abstract} We prove that knowing
the lengths of geodesics
joining points of the boundary of a two-dimensional, compact, simple Riemannian manifold with boundary, we can determine uniquely the Riemannian metric up
to the natural obstruction.
\end{abstract}

\section{Introduction and statement of the results}  \label{section1}
\setcounter {equation}{0}

Let $(M,g)$ be a compact Riemannian manifold with
boundary $\partial M$.
Let $d_g(x,y)$ denote the geodesic distance between $x$ and $y$.
The inverse problem we address in this paper is whether we can determine
the Riemannian metric $g$ knowing $d_g(x,y)$ for any $x\in\partial M$, $y\in
\partial M$.  This problem
arose in rigidity questions in Riemannian geometry \cite{M},
\cite{C}, \cite{Gr}. For the case in which $M$ is a bounded domain
of Euclidean space and the metric is conformal to the Euclidean
one, this problem is known as the inverse kinematic problem which
arose in Geophysics and has a long history (see for instance
\cite{R} and the references cited there).

The metric $g$ cannot be determined from this
information alone. We have $d_{\psi^*g}=d_g$ for any diffeomorphism $\psi:
M\to M$ that leaves the boundary pointwise fixed, i.e.,
$\psi|_{\partial M}=\Id$, where $\Id$ denotes the identity map and $\psi^*g$
is the pull-back of the metric $g$.
The
natural question is whether this is the only obstruction to
unique identifiability of the metric.
It is easy to see that this is not the case.  Namely one can construct a
metric $g$
and find a point
$x_0$
in
$M$
so that
$d_g(x_0, \partial M)> \hbox{ sup }_{x,y \in \partial M} d_g(x,y)$.
For such a metric,
$ d_g $
is independent of a change of
$ g $
in a neighborhood of
$ x_0 $.  The hemisphere of the round sphere is another example.

Therefore it is necessary to impose some a-priori restrictions
on the metric. One such restriction is to assume that the
Riemannian manifold is {\bf simple}, i.e., given two points there is
a unique geodesic joining the points and
$\partial M$
is strictly convex.
$\partial M$
is strictly convex if the second fundamental form of the boundary is
positive definite
in every boundary point.

R.~Michel conjectured in \cite{M} that simple manifolds are
boundary distance rigid that is $d_g$ determines $g$ uniquely up
to an isometry which is the identity on the boundary. This is
known for simple subspaces of Euclidean space (see \cite{Gr}),
simple subspaces of an open hemisphere in two dimensions (see
\cite{M}), simple subspaces of symmetric spaces of constant negative
curvature [BCG], simple two dimensional spaces of negative curvature
(see \cite{C1} or \cite{O}).

In this paper we prove that simple two dimensional compact
Riemannian manifolds are boundary distance rigid. More precisely
we show

\begin{theorem}
Let $(M,g_i), i=1,2$ be a two dimensional simple Riemannian
compact Riemannian manifold. Assume
$$d_{g_1}(x,y)=d_{g_2}(x,y)
\quad   \forall (x,y)\in\partial M\times\partial M
$$
then there exists a diffeomorphism $\psi:  M\to M$,
$\psi|_{\partial M}=Id$, so that
$$g_2=\psi^*g_1.$$
\end{theorem}

As it has been shown in \cite{Sh}, Theorem 1.1 follows from

\begin{theorem}
Let $(M,g_i), i=1,2$ be a two dimensional simple Riemannian
compact Riemannian manifold. Assume
$$d_{g_1}(x,y)=d_{g_2}(x,y) \quad   \forall (x,y)\in\partial
M\times\partial M$$ and
$$g_1|_{\partial M}=g_2|_{\partial M}.$$
then there exists a diffeomorphism $\psi:  M\to M$,
$\psi|_{\partial M}=Id$, $\psi'|_{\partial M}=Id$, so that
$$g_2=\psi^*g_1.$$
\end{theorem}

We will prove Theorem 1.2. The function $d_g$ measures
the travel times of geodesics joining points of the boundary. In
the case that both $g_1$ and $g_2$ are conformal to the Euclidean
metric $e$ (i.e., $(g_k)_{ij}= \alpha_k\delta_{ij}$, $k=1,2$ with
$\delta_{ij}$ the Kr\"onecker symbol), as mentioned earlier,  the
problem we are considering here is known in seismology as the
inverse kinematic problem. In this case, it has been proven by
Mukhometov in two dimensions \cite{Mu} that if $(M,g_i), i=1,2$ is
simple and $d_{g_1}=d_{g_2}$, then $g_1=g_2$. More generally the
same method of proof shows that if $(M, g_i), i=1,2$ are simple
compact Riemannian manifolds with boundary and they are in the
same conformal class, i.e. $g_1=\alpha g_2$ for a positive
function $\alpha$ and $d_{g_1}=d_{g_2}$ then $g_1=g_2$ \cite{Mu1}.
In this case the diffeomorphism $\psi$ must be
the identity. For related results and generalizations see
\cite{B}, \cite{BG}, \cite{C}, \cite{GN}, \cite{MR}.

We mention a closely related inverse problems.
Suppose we have a Riemannian metric in Euclidean space which is the Euclidean metric
outside a compact set. The inverse scattering problem for metrics is
to determine the Riemannian metric by measuring the scattering operator
(see \cite{G}). A similar obstruction occurs in this
case with $\psi$ equal to the identity outside a compact set. It was
proven in \cite{G} that from the wave front set of the scattering operator one can determine,
under some non-trapping assumptions on the metric, the {\bf scattering
relation} on the boundary
of a large ball. We proceed to define in more detail the scattering relation
and its relation with the boundary distance function.

Let $\nu$ denote the unit-inner normal to $\partial M.$ We denote
by $\Omega \left( M\right) \rightarrow M$  the unit-sphere bundle
over $M$:
$$\Omega(M) =\bigcup\limits_{x\in
M}\Omega _{x},\quad \Omega _{x}=\{\xi \in T_{x}(M):\left| \xi
\right|_g =1\}.$$ $\Omega(M)$ is a $(2\mbox{ dim }M-1)$-dimensional compact
manifold with boundary, which can be written as the union
$\partial \Omega \left( M\right) =\partial _{+}\Omega \left(
M\right) \cup
\partial _{-}\Omega \left( M\right) $
$$\partial _{\pm }\Omega \left( M\right)
=\{(x,\xi )\in \partial \Omega \left( M\right) ,\;\pm \,(\nu
\left( x\right) ,\xi )\geq 0\;\}.$$
The manifold of inner vectors
$\partial _{+}\Omega \left( M\right) $ and outer vectors $\partial
_{-}\Omega \left( M\right) $ intersect at the set of tangent
vectors
$$\partial _{0}\Omega \left( M\right)
=\{(x,\xi )\in \partial \Omega \left( M\right),\quad (\nu \left(
x\right) ,\xi )=0\;\}.$$

Let $(M,g)$ be an n-dimensional compact manifold with boundary. We
say that $(M,g)$ is {\bf non-trapping} if each maximal geodesic is
finite. Let $(M,g)$ be non-trapping and the boundary $\partial M$
is strictly convex. Denote by $\tau(x,\xi)$ the length of the
geodesic $\gamma(x,\xi,t), t\geq 0$, starting at the point $x$ in
the direction $\xi \in \Omega_{x}$. These function is smooth on
$\Omega(M)\setminus
\partial_{0}\Omega(M)$. The function $\tau^{0}=\tau|_{\partial\Omega(M)}$
is equal zero on $\partial_{-}\Omega(M)$ and is smooth on
$\partial_{+}\Omega(M)$. Its odd part with respect to $\xi$
$$\tau _{-}^{0}(x,\xi )=\frac{1}{2}\left( \tau
^{0}(x,\xi )-\tau ^{0}\left( x,-\xi \right) \right)$$ is the
smooth function.
\begin{definition}
Let $(M,g)$ be non-trapping with strictly convex boundary. The
scattering relation $\alpha :\partial \Omega \left( M\right)
\rightarrow
\partial \Omega \left( M\right)$ is defined by
$$\alpha (x,\xi )=
(\gamma (x,\xi ,2\tau _{-}^{0}(x,\xi )),\dot{\gamma}(x,\xi ,2\tau
_{-}^{0}(x,\xi ))).
$$
\end{definition}

The scattering relation is a diffeomorphism $\partial \Omega
\left( M\right) \rightarrow \partial \Omega \left( M\right).$
Notice that $\alpha |_{\partial _{+}\Omega \left( M\right)
}:\partial _{+}\Omega \left( M\right) \rightarrow
\partial _{-}\Omega \left( M\right) ,$ $\alpha |_{\partial
_{-}\Omega \left( M\right) }:\partial _{-}\Omega \left( M\right)
\rightarrow \partial _{+}\Omega \left( M\right) $ are
diffeomorphisms as well. Obviously, $\alpha $ is an involution,
$\alpha ^{2}=id$ and $\partial _{0}\Omega \left( M\right) $ is the
hypersurface of its fixed\ points,\ $\alpha (x,\xi )=(x,\xi
),\;(x,\xi )\in \partial _{0}\Omega \left( M\right) .$

A natural inverse problem is whether the scattering relation determines
the metric $g$ up to an isometry which is the identity on the boundary.
In the case that $(M,g)$ is a simple manifold, and we know the metric at the boundary, knowing the scattering
relation is equivalent to knowing the boundary distance function
(\cite{M}). We show in this paper that
if we know the scattering relation we can determine the Dirichlet-to-Neumann (DN) map
associated to the Laplace-Beltrami operator of the metric.
We proceed to define the DN map.

Let $(M,g)$ be a  compact Riemannian manifold with boundary. The
Laplace-Beltrami operator associated to the metric $g$ is given in local
coordinates by
$$
\Delta_ g u=\frac{1}{\sqrt{\det
g}}\sum^n_{i,j=1}\frac{\partial}{\partial x_i}\left(\sqrt{\det g}
g^{ij} \frac{\partial u}{\partial x_j}\right)$$ where $(g^{ij})$
is the inverse of the metric $g$. Let us consider the Dirichlet
problem
$$ \Delta_ g u  = 0\hbox{ on }M,\quad u\Big|_{\partial M}
=  f.$$
We define the DN map in this case by
$$\Lambda_ g(f)=(\nu,\nabla u|_{\partial M})$$
The inverse problem is to recover $g$ from $\Lambda_g.$

In the two dimensional case
the Laplace-Beltrami operator is conformally invariant. More precisely
\[
\Delta_{\beta g}=\frac{1}{\beta}\Delta_g
\]
for any function $\beta$, $\beta\ne 0$. Therefore we have that for $n=2$
\[
\Lambda_{\beta(\psi^\ast g)}=\Lambda_{g}
\]
for any non-zero $\beta$ satisfying $\beta|_{\partial M}=1.$

Therefore the best that one can do in two dimensions is to show that
we can determine the conformal class of the metric $g$ up to
an isometry which is the identity on the boundary
That this is the case is a result proven in [LeU] for simple
metrics and for general connected two dimensional Riemannian manifolds with boundary in [LaU].

In this paper we prove:

\begin{theorem}
Let $(M,g_i), i=1,2$ be compact, simple two dimensional Riemannian
manifolds with boundary.  Assume that $\alpha_{g_1}=\alpha_{g_2}.$ Then
$\Lambda_{g_1}= \Lambda_{g_2}$.
\end{theorem}

The proof of Theorem 1.2 is reduced then to the proof of Theorem
1.3. In fact from Theorem 1.3 and the result of [LaU] we get that
we can determine the conformal class of the metric up to an
isometry which is the identity on the boundary. Now by
Mukhometov's result we have that the conformal factor must be one
proving that the metrics are isometric via a diffeomorphism which
is the identity at the boundary. In other words $d_{g_1}=d_{g_2}$
implies that $\alpha_{g_1}=\alpha_{g_2}$. By Theorem 1.3
$\Lambda_{g_1}=\Lambda_{g_2}.$ By the result of [LeU], [LaU],
there exists a diffeomorphism $\psi:  M\longrightarrow M$,
$\psi|_{\partial M}=\mbox{Identity}$ and a function $\beta\ne 0,
\beta|_{\partial M}=\mbox{identity}$ such that $g_1= \beta
\psi^{\ast} g_2.$ By Mukhometov's theorem $\beta=1$ showing that
$g_1=\psi^{\ast} g_2$ proving Theorem 1.2. and Theorem 1.1.

The proof of Theorem 1.3 consists in showing that from the
scattering relation we can determine the traces at the boundary of
conjugate harmonic functions, which is equivalent information to
knowing the DN map associated to the Laplace-Beltrami operator.
The steps to accomplish this are outlined below. It relies in a
connection between the Hilbert transform and geodesic flow.

We embed $(M,g)$ into a compact Riemannian manifold $(S,g)$ with
no boundary.  Let $\varphi_{t}$ be the geodesic flow on
$\Omega(S)$ and $\mathcal{H}=\frac{d}{dt}\varphi_{t}|_{t=0}$ be
the geodesic vector field. Introduce the map
$\psi:\Omega(M)\rightarrow
\partial_{-}\Omega(M)$ defined by
\[
\psi(x,\xi)=\varphi_{\tau(x,\xi)}(x,\xi), \quad (x,\xi)\in
\Omega(M).
\]
The solution of the boundary value problem for the transport
equation
$$\mathcal{H}u=0,\quad u|_{\partial _{+}\Omega (M)}=w$$
can be written in the form $$u=w_{\psi}=w\circ\alpha\circ \psi.$$

Let $u^{f}$ be the solution of the boundary value problem
\[
\mathcal{H}u=-f,\quad u|_{\partial _{-}\Omega(M)}=0,
\]
 which we can write as
\[
u^{f}(x,\xi)=\int\limits_{0}^{\tau(x,\xi)}f(\varphi_{t}(x,\xi))dt,\quad (x,\xi)\in \Omega(M).
\]
In particular
\[
\mathcal{H}\tau=-1.
\]
The trace
\[
If=u^{f}|_{\partial _{+}\Omega(M)}
\]
is called {\bf the geodesic X-ray transform} of the function $f$. By the fundamental
theorem of calculus we have
\begin{equation}
I\mathcal{H}f=(f\circ\alpha-f) |_{\partial _{+} \Omega(M)}.
\end{equation}
 In what follows we will consider the
operator $I$ acting only on functions that do not depend on $\xi$, unless otherwise
indicated. Let
$L^{2}_{\mu}(\partial_{+}\Omega(M))$ is the real Hilbert space,
with scalar  product given by
\[
(u,v)_{L^{2}_{\mu}(\partial_{+}\Omega(M))}=\int_{\partial_{+}\Omega(M)}\mu uvd\Sigma,\quad \mu=(\xi,\nu).
\]
Here the measure $d\Sigma= d(\partial M)\wedge d\Omega_x$ where
$d(\partial M)$ is the induced volume form on the boundary by the standard measure
on $M$ and
$$d\Omega_x=\sum_{k=1}^{n} (-1)^{k+1} \xi^k d\xi^1\wedge...\wedge \hat{d\xi^k}\wedge...d\xi^n.$$

As usual the scalar product in $L^{2}(M)$
is defined by $$(u,v)=\int_{M}uv\sqrt{detg}dx.$$  The operator $I$ is a
bounded operator from $L^{2}(M)$ into
$L^{2}_{\mu}(\partial_{+}\Omega(M))$. The adjoint
$I^{\ast}:L^{2}_{\mu}(\partial_{+}\Omega(M))\rightarrow L^{2}(M)$
is given by
\[
I^{\ast}w(x)=\int_{\Omega_{x}}w_{\psi}(x,\xi)d\Omega_{x},
\]

We will study the solvability of equation $I^{\ast}w=h$ with
smooth right hand side. Let $w\in
C^{\infty}(\partial_{+}\Omega(M))$. Then the function $w_{\psi}$
will not be smooth on $\Omega(M)$ in general. We have that
$w_{\psi}\in C^{\infty}
(\Omega(M)\setminus\partial_{0}\Omega(M))$. We give below
necessary and sufficient conditions for smoothness of $w_{\psi}$
on $\Omega(M)$.

We introduce the operators of even and odd continuation
with respect to $\alpha $:
\[
A_{\pm }w(x,\xi )=w(x,\xi ),\;\;\;\;\;\;
(x,\xi )\in \partial _{+}\Omega \left( M\right),
\]
\[
A_{\pm }w(x,\xi )=\pm \left( \alpha ^{\ast }w\right) (x,\xi ),
\;\;(x,\xi )\in \partial _{-}\Omega \left( M\right).
\]

The scattering relation preserves the measure $|(\xi, \nu)|
d\Sigma$ and therefore the operators
$A_{\pm}:L^2_{\mu}(\partial_{+}\Omega(M))\rightarrow L_{\left| \mu
\right| }^{2}\left( \partial\Omega \left( M\right) \right)$ are
bounded, where $L_{\left| \mu \right| }^{2}\left(\partial \Omega
\left( M\right) \right) $ is real Hilbert space with scalar
product
\[
\left( u,v\right) _{L_{\left| \mu \right| }^{2}\left(\partial
\Omega \left( M\right) \right) }=\int\limits_{\partial \Omega
\left( M\right) }\left| \mu \right| uvd\Sigma,\quad \mu=(\xi,\nu).
\]
The adjoint of $A_{\pm}$ is a bounded operator
$A^*_{\pm}:L_{\left| \mu \right| }^{2}\left(\partial \Omega \left(
M\right) \right)\rightarrow L^2_{\mu}(\partial_{+}\Omega (M))$
given by
\[
A_{\pm}^{\ast}u=(u\pm u\circ \alpha)|_{\partial_{+}\Omega(M)}.
\]
Using $A_{-}^{\ast}$ formula (1.1) can be written in the form
\begin{equation}
I\mathcal{H}f=-A_{-}^{\ast}f^{0}, \quad f^{0}=f|_{\partial \Omega
(M)}.
\end{equation}

 The space $C_{\alpha }^{\infty }\left(
\partial _{+}\Omega \left( M\right) \right) $ is defined by
\[
C_{\alpha }^{\infty }\left( \partial _{+}\Omega \left( M\right) \right) =\{w\in C^{\infty }
\left( \partial _{+}\Omega \left( M\right) \right)
 :w_{\psi }\in C^{\infty }\left( \Omega \left( M\right) \right) \}.
\]
We have the following characterization of the space of smooth solutions
of the transport equation
\begin{lemma}
$$C_{\alpha}^{\infty}(\partial_{+}\Omega(M))=
\{w\in
C^{\infty}(\partial_{+}\Omega(M)):A_{+}w\in
C^{\infty}(\partial\Omega(M))\}.$$
\end{lemma}

Now we can state the main theorem for solvability for $I^{\ast}$.
\begin{theorem}
Let $(M,g)$ be a simple, compact two dimensional Riemannian manifold with boundary. Then the
operator
$I^{\ast}:C_{\alpha}^{\infty}(\partial_{+}\Omega(M))\rightarrow
C^{\infty}(M)$ is onto.
\end{theorem}

Now we define the Hilbert transform:
\begin{equation}
Hu(x,\xi)=\frac{1}{2\pi}\int_{\Omega_{x}}\frac{1+(\xi,\eta)}{(\xi_{\perp},\eta)}u(x,\eta)d\Omega_{x}(\eta),\quad
\xi\in \Omega_{x},
\end{equation}
where the integral is understood as a principle value integral.
Here $\perp $\ means\ a $90^o$ degree rotation.\ In coordinates
$(\xi _{\perp })_{i}=\varepsilon _{ij}\xi ^{j},$\
\vspace{1pt}where
\[
\varepsilon =\sqrt{\det g}\begin{pmatrix}
  _{0} & _{1} \\
  _{-1} & _{0}
\end{pmatrix}.
\]

The Hilbert transform $H$ transforms even (respectively odd) functions with respect
to $\xi$ to even (respectively  odd) ones. If $H_{+}$ (respectively
$H_{-})$ is the even (respectively odd) part of
the operator $H$:
\[
H_{+}u(x,\xi)=\frac{1}{2\pi}\int_{\Omega_{x}}\frac{(\xi,\eta)}{(\xi_{\perp},\eta)}u(x,\eta)d\Omega_{x}(\eta),
\]
\[
Hu_{-}(x,\xi)=\frac{1}{2\pi}\int_{\Omega_{x}}\frac{1}{(\xi_{\perp},\eta)}u(x,\eta)d\Omega_{x}(\eta)
\]
and $u_{+}, u_{-}$ are the even and odd parts of the function $u$, then
$H_{+}u=Hu_{+}, H_{-}u=Hu_{-}$.
\smallskip

We introduce the notation
$\mathcal{H}_{\perp}=(\xi_{\perp},\nabla)=-(\xi,\nabla_{\perp})$,
where $\nabla_{\perp}=\varepsilon\nabla$ and $\nabla$ is the
covariant derivative with respect to the metric $g$. The following
commutator formula for the geodesic vector field  and the Hilbert
transform is very important in our approach.

\begin{theorem}
Let $(M,g)$ be a two dimensional Riemannian manifold.
For any smooth function $u$ on $\Omega(M)$ we have
the identity
\begin{equation}
[H,\mathcal{H}]u=\mathcal{H}_{\perp}u_{0}+(\mathcal{H}_{\perp}u)_{0}
\end{equation}
where
\[
u_{0}(x)=\frac{1}{2\pi}\int_{\Omega_{x}}u(x,\xi)d\Omega_{x}
\]
is the average value.
\end{theorem}

\smallskip

Now we can prove the Theorem 1.3.

\bigskip

Separating the odd and even parts with respect to $\xi$ in (1.4)
we obtain the identities:
\[
H_{+}\mathcal{H}u-\mathcal{H}H_{-}u=(\mathcal{H}_{\perp}u)_{0},\quad
H_{-}\mathcal{H}u-\mathcal{H}H_{+}u={\mathcal{H}}_{\perp}u_{0}.
\]

Let $(M,g)$ be a non-trapping strictly convex manifold. Take
$u=w_{\psi},w\in C_{\alpha}^{\infty}(\partial_{+}(\Omega))$. Then
$$2\pi \mathcal{H}H_{+}w_{\psi}=-\mathcal{H}_{\perp}I^{\ast}w$$ and using
the formula (1.2) we conclude
\begin{equation}
2\pi A_{-}^{\ast}H_{+}A_{+}w=I\mathcal{H}_{\perp}I^{\ast}w,
\end{equation}
since $w_{\psi }|_{\partial \Omega (M)}=A_{+}w$.

Let $(h,h_{\ast})$ be a pair of conjugate harmonic functions on
$M$,
\[
\nabla h=\nabla_{\perp}h_{\ast},\quad \nabla
h_{\ast}=-\nabla_{\perp}h.
\]

Notice, that $\delta\nabla=\triangle$ is the Laplace-Beltrami
operator and $\delta\nabla_{\perp}=0$. Let $I^{\ast}w=h$. Since
$I\mathcal{H}_\perp h=I\mathcal{H}
h_{\ast}=-A^{\ast}_{-}h^{0}_{\ast}$, where
$h^{0}_{\ast}=h_{\ast}|_{\partial M}$, we obtain from (1.5)
\begin{equation}
2\pi A_{-}^{\ast}H_{+}A_{+}w=-A^{\ast}_{-}h_{\ast}^{0}.
\end{equation}

The following theorem gives the key to obtain the DN map from the
scattering relation.

\smallskip
\begin{theorem}
Let $M$ be a 2-dimensional simple manifold. Let $w\in
C^{\infty}_{\alpha}(\partial_{+}\Omega(M))$ and $h_{\ast}$ is
harmonic continuation of function $h_{\ast}^{0}$. Then the
equation (1.6) holds iff the functions $h=I^{\ast}w$ and $h_{\ast}$ are
conjugate harmonic functions.
\end{theorem}

\begin{proof}

The necessity has already been established. Using (1.2) and (1.5)
the equality (1.6) can be written in the form
\[
I\mathcal{H}_{\perp}h=I\mathcal{H}q,
\]
where $q$ is an arbitrary smooth continuation onto $M$ of the
function $h_{\ast}^{0}$ and $h=I^{\ast}w$. Thus, the ray transform
of the vector field $\nabla q+\nabla _{\perp }h$ equal $0$.
Consequently, this field is potential ([An]), that is, $\nabla
q+\nabla _{\perp }h=\nabla p$ and $\ p|_{\partial M}=0.$ Then the
functions $h$ and $h_{\ast}=q-p$ are conjugate harmonic functions
and $h_{\ast}|_{\partial M}=h^{0}_{\ast}$. We have finished the
proof of the main theorem.

\end{proof}

\smallskip

In summary we have the following procedure to obtain the DN map
from the scattering relation. For an arbitrary given smooth
function $h_{\ast}^{0}$ on $\partial M$ we find a solution $w\in
C^{\infty}_{\alpha}(\partial_{+}\Omega(M))$ of the equation (1.6).
Then the functions $h^{0}=2\pi (A_{+}w)_{0}$ (notice, that $2\pi
(A_{+}w)_{0}=I^{\ast}w|_{\partial M} $) and $h_{\ast}^{0}$ are the
traces of conjugate harmonic functions. It gives the map
\[
h_{\ast}^{0}\rightarrow
(\nu_{\perp},\nabla^{\parallel}h^{0})=(\nu,\nabla
h_{\ast})|_{\partial M}.
\]
which is the DN map. This proves Theorem 1.3. For the definition of
$\nabla^{\parallel}$ see section 2.

A brief outline of the paper is as follows. In section 2 we
collect some facts and definition we need later. In section 3 we
study the solvability of $I^{\ast}w=h$ on Sobolev spaces and prove
Theorem 1.4. In section 4 we make a detailed study of the
scattering relation and prove Lemma 1.1. In section 5 we prove
Theorem 1.5.


\section{\bf Preliminaries and Notation}
\setcounter {equation}{0}

Here we will give some definitions and formulas, that we need in
what follows. For further references see \cite{E}, \cite{J}, \cite{K},
\cite{Sh}. Let $\pi:T(M)\rightarrow M$ be the tangent bundle over
an $n$-dimensional Riemannian manifold $(M,g)$. We will denote points
of the manifold $T(M)$ by pairs $(x,\xi)$. The connection map
$K:T(T(M))\rightarrow T(M)$ is defined by its local representation
$$K(x,\xi,y,\eta)=(x,\eta+\Gamma(x)(y,\xi)),\quad(\Gamma(x)(y,\xi))^{i}=\Gamma^{i}_{jk}(x)y^{j}\xi^{k}, i=1,...,n,$$
where $\Gamma^{i}_{jk}$ are the Christoffel symbols of the metric $g$

$$\Gamma^{i}_{jk}=\frac{1}{2}g^{il}(\frac{\partial g_{jl}}{\partial x^{k}}+\frac{\partial g_{kl}}{\partial x^{j}}
-\frac{\partial g_{jk}}{\partial x^{l}}).$$ The linear map
$K(x,\xi)=K|_{(x,\xi)}:T_{(x,\xi)}(T(M))\rightarrow T_{x}M$
defines the horizontal subspace $H_{(x,\xi)}=KerK(x,\xi)$. It can
be identified with the tangent space $T_{x}(M)$ using the operator
$$J^{h}_{(x,\xi)}=(\pi'(x,\xi)|_{H_{(x,\xi)}})^{-1}:T_{x}(M)\rightarrow
H_{(x,\xi)}.$$ The vertical space $V_{(x,\xi)}=Ker\pi'(x,\xi)$ can
also identified with $T_{x}(M)$ by using the operator
$$J^{v}_{(x,\xi)}=(K(x,\xi)|_{V_{(x,\xi)}})^{-1}:T_{x}(M)\rightarrow V_{(x,\xi)}.$$
The tangent space $T_{(x,\xi)}(T(M))$ is the direct sum of
the horizontal and vertical subspaces,
$T_{(x,\xi)}(T(M))=H_{(x,\xi)}\oplus V_{(x,\xi)}$. An arbitrary
vector $X\in T _{(x,\xi)}(T(M))$ can be uniquely decomposed in the form
$$X=J^{h}_{(x,\xi)}X_{h}+J^{v}_{(x,\xi)}X_{v},$$
where
$$X_{h}=\pi'(x,\xi)X,\quad X_{v}=K(x,\xi)X.$$
We will call $X_{h},X_{v}$ the horizontal and vertical components
of the vector $X$ and use the notation $X=(X_{h},X_{v})$. If in
local coordinates $X=(X^{1},...,X^{2n})$ then $X_{h},X_{v}$ is
given by
$$X_{h}^{i}=X^{i},\quad X_{v}^{i}=X^{i+n}+\Gamma^{i}_{jk}(x)X^{j}\xi^{k},\quad i=1,...,n.$$
The horizontal and vertical components $Y^{h},Y^{v}\in
T^{\ast}_{x}(M)$ of a covector $Y\in T^{\ast} _{(x,\xi)}(T(M)),
Y=(Y^{h},Y^{v})$ are given by
$$\langle X,Y\rangle=(X_{h},Y^{h})+(X_{v},Y^{v}).$$
In local coordinates we have
$$Y^{h}_{j}=Y_{j}-\Gamma^{i}_{jk}(x)Y_{i+n}\xi^{k},\quad Y^{v}_{j}=Y_{j+n},\quad j=1,...,n.$$

\smallskip

Let $N$ be a smooth manifold and $f:T(M)\rightarrow N$ a smooth
map. Then the derivative $f'(x,\xi):T_{(x,\xi)}(T(M))\rightarrow
T_{f(x,\xi)}(N)$ defines the horizontal $\nabla_{h} f(x,\xi)$ and
vertical $\nabla_{v} f(x,\xi)$ derivatives:
\[
\nabla_{h} f(x,\xi)=f'(x,\xi)\circ
J^{h}_{(x,\xi)}:T_{x}(M)\rightarrow T_{f(x,\xi)}(N),
\]
\[
\nabla_{v} f(x,\xi)=f'(x,\xi)\circ
J^{v}_{(x,\xi)}:T_{x}(M)\rightarrow T_{f(x,\xi)}(N).
\]
We have then that
\begin{equation}
f'(x,\xi)X=(\nabla_{h}f(x,\xi),X_{h})+(\nabla_{v}f(x,\xi),X_{v})
\end{equation}
In local coordinates
\[
\nabla_{h j} f^{(\alpha)}(x,\xi)=(\frac{\partial}{\partial x^{j}}-
\Gamma^{i}_{jk}(x)\xi^{k}\frac{\partial}{\partial
\xi^{i}})f^{(\alpha)}(x,\xi),
\]
\[
\nabla_{v j} f^{(\alpha)}(x,\xi)=\frac{\partial}{\partial
\xi^{j}}f^{(\alpha)}(x,\xi),\quad\alpha=1,...,\mbox{dim}N.
\]

\bigskip
Let $e:\Sigma\rightarrow T(M)$ be a smooth embedding and
$p:T(M)\rightarrow \Sigma$ denotes a normal projection, i.e. $p
e=id.$ For a point $(x,\xi)\in e(\Sigma)$, $\mbox{ Ker
}p'(x,\xi)=N(\Sigma)\bigcap T_{(x,\xi)}(T(M))$, where $N(\Sigma)$
is the normal bundle of $\Sigma$. 

\begin{definition}
The horizontal and vertical derivatives $\nabla_{h}^{\parallel},
\nabla_{v}^{\parallel}$ on $\Sigma$ are the tangent components of
$\nabla_{h}$ and $\nabla_{v}$ with respect to $\Sigma$:
\[
\nabla^{\parallel}_{h}f=\nabla_{h}(f\circ p)|_{\Sigma},\quad
\nabla_{v}^{\parallel}f=\nabla_{v}(f\circ p)|_{\Sigma}.
\]
\end{definition}

This definition is obviously independent of the choice of the
normal projection. From the identity $f=(f\circ p)\circ e$ and
(2.1) we have for any vector $X\in T_{(x,\xi)}(T(\Sigma)),$ by
identifying $X$ with its embedding $e'(x,\xi)X$, that
\[
f'(x,\xi)X=(f\circ p)'(x,\xi)X=(\nabla_{h}(f\circ
p)(x,\xi),X_{h})+(\nabla_{v}(f\circ p)(x,\xi),X_{v}).
\]
Therefore
\[
f'(x,\xi)X=(\nabla^{\parallel}_{h}f(x,\xi),X_{h})+(\nabla^{\parallel}_{v}f(x,\xi),X_{v}).
\]
Let us give an equivalent definition. Denote by $i( x,\xi)$
the isomorphism
\[
i\left( x,\xi \right) :T_{\left( x,\xi \right)
}\left( T\left( M\right) \right) \rightarrow T_{x}(M)\times
T_{x}\left( M\right),
\mbox{ defined by }
i\left( x,\xi \right) X=(X_{h},X_{v})
\]
Let $i_{\Sigma }\left( x,\xi \right) $ be its restriction on
$T_{\left( x,\xi \right) }\left( \Sigma \right) ,\;i_{\Sigma
}\left( x,\xi \right) =i\left( x,\xi \right) |_{T_{\left( x,\xi
\right) }\left( \Sigma \right) }$. \ Denote by $T_{x}^{\left(
i\right) }\left( M\right) =pr_{i}\;Ran\,i_{\Sigma }\left( x,\xi
\right) \subset T_{x}(M)\times T_{x}\left( M\right) ,\;i=1,2,$\
the projections to the first and second component of the range of
$i_\Sigma.$ Then
\[\nabla
_{h}^{\parallel }f\left( x,\xi \right) =\nabla _{h}\tilde{f}\left(
x,\xi \right) |_{T_{x}^{\left( 1\right) }\left( M\right)
},\;\;\;\nabla _{v}^{\parallel }f\left( x,\xi \right) =\nabla
_{v}\tilde{f}\left( x,\xi \right) |_{T_{x}^{\left( 2\right)
}\left( M\right) },
\]
where $\tilde{f}:T\left( M\right) \rightarrow N$\ is an arbitrary
smooth continuation of $f$.

\smallskip
\noindent{\bf Example 1}

Let us consider the special case of $\Sigma=\Omega(M)$ the unit
sphere bundle. Then $p(x,\xi)=(x,\xi/|\xi|)$ is a normal
projection. For simplicity of notation in the case of the manifold
$\Omega(M)$ we will use notation $\partial$ instead of
$\nabla_{v}^{\parallel}$. Then
$$\partial=(\frac{\partial}{\partial\xi}-\xi(\xi,\frac{\partial}{\partial\xi}))|_{|\xi|=1}.$$
Clearly, $(\xi,\partial)=0$. Notice, that the horizontal
derivative is tangent to the submanifold $\Omega(M)$ since
$\nabla_{h}|\xi|=0$. Because of this we keep the original notation
for horizontal derivative on $\Omega(M)$ and moreover simplify
this notation to just $\nabla$. Thus, in the case of the manifold
$\Omega(M)$ we will use notation $\nabla$ and $\partial$ for the
horizontal and vertical derivatives. We will also use $\nabla$ for
horizontal derivative on $T(M)$. Notice relation
\[
T_{(x,\xi)}(\Omega(M))=\{X\in T_{(x,\xi)}(T(M)):(\xi,X_{v})=0\}.
\]

\smallskip
\noindent{\bf Example 2}

Let $\Gamma$ be a smooth hypersurface in $M$ with normal $\nu$ and
$\Sigma$ its lift to $\Omega(M)$,
$\Sigma=\pi^{-1}_{1}(\Gamma)=\{(x,\xi):x\in \Gamma,
\xi\in\Omega_{x}(M)\}$, where $\pi_{1}=\pi|_{\Omega(M)}$. Then
\begin{equation}
T_{(x,\xi)}(\pi^{-1}_{1}(\Gamma))=\{X\in
T_{(x,\xi)}(T(M)):(\nu(x),X_{h})=0,(\xi,X_{v})=0\}
\end{equation}
It means that $T_{x}^{\left( 1\right) }\left( M\right) =\{\eta \in
T_{x}\left( M\right) :\left( \eta ,\nu \left( x\right) \right)
=0\},\;T_{x}^{\left( 2\right) }(M)=\{\eta \in T_{x}\left( M\right)
:\left( \xi ,\eta \right) =0\}$. We have the same vertical
derivative as in the case of $\Omega \left( M\right) .$\ The
horizontal derivative for this case will be denoted by $\nabla
^{\parallel },$ so that
\[
\nabla^{\parallel}=(\nabla-\nu(\nu,\nabla))|_{\pi^{-1}_{1}(\Gamma)}.
\]
Clearly, $(\nu,\nabla^{\parallel})=0$.

We now state a similar definition of vertical and horizontal
derivatives on a submanifold for semibasic tensor fields. We
recommend chapter 3 of [Sh] for more details.

Let $T_{s}^{r}(M)$ denote the  bundle of tensor fields of degree
$(r,s)$ on $M.$  A section of this bundle is called a tensor field
of degree $(r,s)$.
Let $\pi _{s}^{r}:T_{s}^{r}(M)\rightarrow M$ be the projection.
A\ fiber map $u:T(M)\rightarrow
T_{s}^{r}(TM)$, i.e. $\pi _{s}^{r}\circ u=\pi $ is called a semibasic
tensor field of degree $(r,s)$ on the manifold $T(M)$. Denote by
$\xi $ the semibasic vector field given by the identity map
$T(M)\rightarrow T(M)$. An arbitrary tensor field $u$ of degree $(r,s)$
on the manifold $M$, i.e. section $u:M\rightarrow T_{s}^{r}(M)$
defines by the formula $u\circ \pi $ a semibasic tensor field (since
$\pi _{s}^{r}\circ (u\circ \pi )=(\pi _{s}^{r}\circ u)\circ \pi
=id\circ \pi =\pi $). The map $u\rightarrow u\circ \pi $
identifies tensor fields on $M$ and $\xi $-constant semibasic
tensor fields on $T(M)$. A semibasic tensor field on the submanifold
$\Sigma \subset T(M)$ is defined as the fiber map $u:\Sigma
\rightarrow T_{s}^{r}(M)$, i.e. $\pi _{s}^{r}\circ u=\pi |_{\Sigma
}$. If $u$ is semibasic tensor field on $T(M)$, then its restriction
$u|_{\Sigma }$ is a semibasic tensor field on $\Sigma $. Using the metric $g$
we can identify the bundle $T_{s}^{r}(M)$ with $T_{0}^{r+s}(M)$
and the bundle $T_{r+s}^{0}(M)$ with $T_{r+s}^{\ast }(M).$

A semibasic tensor field $u:T(M)\rightarrow T_{m}^{0}(TM)$ is\
naturally identified with the polynomial of degree $m$
in $\eta$,  $u_{(m)}:T(M)\times
T^{m}(M)\rightarrow R,$ $u_{m}\left( x,\xi ,\eta _{\left( 1\right)
},...,\eta _{\left( m\right) }\right) =u_{i_{1}...i_{m}}\left(
x,\xi \right) \eta _{\left( 1\right) }^{i_{1}}...\eta _{\left(
m\right) }^{i_{m}},$ where the fibers of the vector bundle
$\ T(M)\times
T^{m}(M)=T^{m+1}(M)$ are given by
$$T_{x}^{m+1}(M)=\underset{m+1}{\text{ }\underbrace{T_{x}\left( M\right)
\times ...\times T_{x}\left( M\right) }}.$$

The horizontal derivative $\nabla _{h}u_{\left( m\right) }$ and
vertical derivative $\nabla _{v}u_{\left( m\right) }$ are the fiber
maps $T(M)\times T^{m}(M)\rightarrow T^{\ast }(M)=T_{1}^{0}(M)$,\
and \ therefore are semibasic tensor fields on $T(M)\times T^{m}(M)$
of  degree 1. The corresponding polynomials $(\nabla _{h}u_{\left(
m\right) })_{\left( 1\right) },(\nabla _{v}u_{\left( m\right)
})_{\left( 1\right) }$ on $T(M)\times T^{m+1}(M)$\ define
semibasic tensor fields $\nabla u,\;\nabla _{\xi }u$\ of degree $(m+1)$,
which we call the horizontal and vertical derivatives of $u.$
Thus, by definition
\[
(\nabla u)_{(m+1)}=
(\nabla _{h}u_{\left( m\right) })_{\left( 1\right) },\;(\nabla
_{\xi }u)_{(m+1)}= (\nabla _{v}u_{\left( m\right) })_{\left(
1\right)}.
\]

Straightforward calculations give the representation in local coordinates:
\[
\left( \nabla u\right) _{i_{1}...i_{m+1}}=\tilde{\nabla}_{i_{m+1}}u_{i_{1}...i_{m}}-
\Gamma _{i_{\left( m+1\right) }k}^{j}\xi ^{k}\frac{\partial u_{i_{1}...i_{m}}}{\partial \xi ^{j}},\quad
(\nabla _{\xi }u)_{i_{1}...i_{m+1}}=
\frac{\partial u_{i_{1}...i_{m}}}{\partial \xi ^{i_{m+1}}},
\]
where $\tilde{\nabla}$ denotes the usual covariant derivative on the
manifold $(M,g).$ Notice, that for $\xi $-constant tensor fields, $\nabla
u=\tilde{\nabla}u$ and since we identify $\xi $-constant semibasic tensor
fields with tensor fields on $M$, we will use one notation $\nabla $\ for
covariant and horizontal derivatives.

As in the case of maps, we define tangent derivatives of semibasic
tensor fields on the submanifold $\Sigma \subset T(M)$:
\[
\nabla ^{\parallel }u=\nabla (u\circ p)|_{\Sigma },\quad \nabla
_{\xi }^{\parallel }u=\nabla _{\xi }(u\circ p)|_{\Sigma }.
\]

In the case $\Sigma =\Omega (M)$ we keep the notation $\nabla $ instead of $\nabla
^{\parallel }$ and use $\partial $ instead $\nabla _{\xi
}^{\parallel }$. We mention the following formulas, [Sh] :

\begin{equation}
\nabla g=0,\quad \nabla \xi =0,\quad \partial _{j}\xi ^{i}= \delta
_{j}^{i}-\xi ^{i}\xi _{j},
\end{equation}
\begin{equation}
\lbrack \nabla ,\partial \rbrack =0,\quad \lbrack \partial
_{i},\partial _{j}\rbrack =\xi _{i}\partial _{j}-\xi _{j}\partial
_{i},
\end{equation}
\begin{equation}
\lbrack \nabla _{i},\nabla _{j}\rbrack u=-R_{qij}^{p}\partial
_{p}u,
\end{equation}
where $R$ is the curvature tensor. In the last formula $u$ is a
scalar. In the case that $\Sigma =\pi _{1}^{-1}\left( \Gamma \right)
$ we will use the same\ notations $\nabla ^{\parallel },\partial $
as in the case of maps.

\section{\bf The geodesic X-ray transform}
\setcounter {equation}{0}

In this section we study the solvability of the equation $I^*w=h$
and prove Theorem 1.4.

\begin{lemma}
Let $V$ be an open set of a Riemannian manifold $(M,g).$
We can define the ray transform as before. Then the
normal operator
$I^{\ast}I$
is an elliptic pseudodifferential operator of order $-1$ on $V$ with
principal symbol $c_{n}\left| \xi \right|
^{-1}$ where $c_{n}$ is a constant.
\end{lemma}

\begin{proof}

It is easy to see, that
\begin{equation}
\left( I^{\ast }If\right) \left( x\right) = \int\limits_{\Omega
_{x}}d\Omega _{x}\int\limits_{-\tau \left( x,-\xi \right) }^{\tau
\left( x,\xi \right) }f\left( \gamma \left( x,\xi ,t\right)
\right) dt= 2\int\limits_{\Omega _{x}}d\Omega
_{x}\int\limits_{0}^{\tau \left( x,\xi \right) }f\left( \gamma
\left( x,\xi ,t\right) \right) dt.
\end{equation}

Before we continue we make a remark concerning notation. We have
used up to know the notation $\gamma (x,\xi ,t)$ for a geodesic.
But it is known \cite{J} , that a geodesic depends smoothly on the
point $x$ and vector $\xi t\in T_{x}(M).$ Therefore in  what
follows we will also use sometimes the notation $\gamma (x,\xi t)$
for a geodesic. Since the manifold  $M$ is simple and any small
enough neighborhood $U$ (in $\left( S,g\right) $) is also simple
(an open domain is simple if its closure is simple). \ For any
point $x\in U$\ \ there is an open domain $D_{x}^{U}\subset
T_{x}\left( U\right) $ such that exponential map
$exp_{x}:D_{x}^{U}\rightarrow U,\;exp_{x}\eta =\gamma (x,\eta )$\
is a diffeomorphism  onto $U.$ Let $D_{x},\;x\in M$ \ be the
inverse image of $M$, then $exp_{x}(D_{x})=M$ and
$exp_{x}|_{D_{x}}:D_{x}\rightarrow M$ is a diffeomorphism.

Now we change variables in (3.1), $y=\gamma (x,\xi t).$ Then
$t=d_g\left( x,y\right)$ and
\[
(I^{\ast }If)\left( x\right) = \int\limits_{M}K\left( x,y\right)
\,f\left( y\right) dy,
\]
where
\[
K\left( x,y\right) =2\frac{\det \left(
exp_{x}^{-1}\right) ^{\prime }\left( x,y\right) \sqrt{\det g\left(
x\right) }}{d_g^{n-1}\left( x,y\right)}.
\]

Notice, that since
\begin{equation}
\gamma (x,\eta )=x+\eta +O(\left| \eta \right| ^{2}),
\end{equation}
it follows, that the Jacobian matrix of the exponential map is $1$
at 0, and then $\det (exp_{x}^{-1\prime }\left)( x,x\right)
=1/\det \left( exp_{x}\right) ^{\prime }\left( x,0\right) =1$.
>From (3.2) we also conclude that
\[
d^{2}\left( x,y\right) =G_{ij}\left( x,y\right) \left( x-y\right) ^{i}\left( x-y\right) ^{j},\;\;G_{ij}\left( x,x\right) =
g_{ij}\left( x\right) ,\;\;G_{ij}\in C^{\infty }\left( M\times M\right)
\]
Therefore the kernel of $I^*I$ can be written in the form
\[
K\left( x,y\right) =\frac{2\det \left( exp_{x}^{-1}\right) ^{\prime }\left( x,y\right)
\sqrt{\det g\left( x\right) }}{\left( G_{ij}\left( x,y\right)
\left( x-y\right) ^{i}\left( x-y\right) ^{j}\right) ^{\left(
n-1\right) /2}}.
\]

Thus \ the kernel $K$ has at the diagonal $x=y$ \ a singularity
of type $\left| x-y\right| ^{-n+1}.$ \ The kernel
\[
K_{0}^{\,}\left( x,y\right) =\frac{2\sqrt{\det g\left( x\right)
}}{\left( g_{ij}\left( x\right) \left( x-y\right) ^{i}\left(
x-y\right) ^{j}\right) ^{\left( n-1\right) /2}}
\]
has the same  singularity. Clearly, the difference $K-K_{0}$ has
a singularity of type $\left| x-y\right| ^{-n+2}.$  Therefore the
principal symbols of both operators coincide.\ The principal
symbol of the integral operator, corresponding to the kernel
$K_{0}$ coincide with its full symbol and is easily\ calculated.\ As
a result

\[
\sigma \left( I^{\ast }I\right) (x,\xi )=
2\sqrt{\det g\left( x\right) }\int \frac{e^{-
i(y,\xi )}}{\left( g_{ij}\left( x\right) y^{i}y^{j}\right) ^{\left( n-1\right) /2}}dy=
c_{n}\left| \xi \right| ^{-1}.
\]

\end{proof}

Let $r_{M}$ denotes the restriction from $S$ onto $M.$

\begin{theorem}
Let $U$ be a simple neighborhood of the simple
manifold $M$. Then for any function $h\in H^{s}\left( M\right)
,\,s\geq 0\;$ there exists function $f\in H^{s-1}\left( U\right)
,\;r_{M}I^{\ast }If=h.$
\end{theorem}

\begin{proof}

Let $(M,g)$ be simple and\ embedded into a compact Riemannian
manifold $\left( S,g\right) \;$ without boundary, of the same
dimension. Choose a\ finite atlas\ of $S$, which consist of simple
open sets $U_{k}$ with coordinate maps $ \kappa
_{k}:U_{k}\rightarrow R^{n}$. Let $\left\{ \varphi _{k}\right\} $
be the subordinated partition of unity:\; $\varphi _{k}\geq
0,\;supp\varphi _{k}\subset U_{k},\;\sum \varphi _{k}=1$.  We
assume without loss of generality that $M\subset U_1$ and $\varphi
_{1}|_{M}=1$. We consider the operators $I_{k},$ $I_{k}^{\ast }$
for the domain $U_{k},$ and the pseudodifferential operator on
$\left( S,g\right) $
\[
Pf=\sum\limits_{k}\varphi _{k}\left( I_{k}I_{k}^{\ast }\right)
\left( f|_{U_{k}}\right) ,\;\;f\in D^{\prime }\left( X\right).
\]
Every operator $I_{k}I_{k}^{\ast }:C_{0}^{\infty }\left(
U_{k}\right) \rightarrow C^{\infty }\left( U_{k}\right) $ is an
elliptic pseudodifferential operator of order $-1$ with principal
symbol $c_{n}\left| \xi \right| ^{-1},\;\xi \in T\left(
U_{k}\right) .$ Then $P$ is an elliptic\ pseudodifferential
operator with principal symbol $c_{n}\left| \xi \right| ^{-1},$
$\xi \in T\left( S\right) $ and, therefore, is a Fredholm operator
from $H^{s}(S) $ into $H^{s+1}(S).$  We have that $\mbox{ Ker }P$
has finite dimension, $\mbox{ Ran }P$ is closed and has finite
codimension.\ Notice, that $P^{\ast }=P$ (more precisely if
$P^{s}=P:H^{s}\left( S\right) \rightarrow H^{s+1}\left( S\right)
,$ then $\left( P_{s}\right) ^{\ast }=P_{-s-1}$)$.$

\smallskip
For arbitrary $s\geq 0$\ the operator $r_{M}:H^{s}\left( S\right)
\rightarrow H^{s}\left( M\right) $ is bounded and
$r_{M}(H^{s}\left( S\right) )=H^{s}\left( M\right) .$ Then the
range of $r_{M}P:H^{s}\left( S\right) \rightarrow H^{s+1}\left(
M\right) ,\;s\geq -1$ is closed.

\smallskip
Since $M$ is only covered by $U_{1}$ and $\varphi _{1}|_{M}=1$ we
have that $r_{M}Pf=r_{M}I_{1}^{\ast }I_{1}\left( f|_{U_{1}}\right)
$. Thus, the\ range of the \ operator$\;r_{M}I_{1}^{\ast
}I_{1}:H^{s}\left( U_{1}\right) \rightarrow H^{s+1}\left( M\right)
,\;s\geq -1$ is closed. Now to prove the\ solvability of the\
equation
\[
r_{M}I_{1}^{\ast }I_{1}f=h\in H^{s+1}\left( M\right) ,\;s\geq -1
\]
in $H^{s}\left( U_{1}\right) $\ is sufficient to show that the kernel of
the adjoint $\left( r_{M}I_{1}^{\ast }I_{1}\right) ^{\ast
}:\left( H^{\left( s+1\right) }\left( M\right) \right) ^{\ast
}\rightarrow \left( H^{s}\left( U_{1}\right) \right) ^{\ast }$ is
zero.

\smallskip
Let $\left\langle ,\right\rangle _{M}$\ and $\left\langle
,\right\rangle $\ be\ dualities between $H^{s}(M)$ and $\left(
H^{s}\right) ^{\ast }(M)$ or $H^{s}(S)$ and $H^{-s}(S)$
respectively. The dual\ space $\left( H^{s}\left( M\right) \right)
^{\ast },\;s\geq 0$ can be identified with the subspace of
$H^{-s}\left( S\right) :$
$$\left( H^{s}\left( M\right) \right) ^{\ast }=H^{-s}\left( M\right) =
\left\{ u\in H^{-s}\left( S\right) :{supp}u\subset M\right\} .$$
For any $f\in H^{s}\left( U_{1}\right) ,\;u\in H^{-(1+s)}\left(
M\right) $ we have
$$\langle r_{M}I_{1}^{\ast }I_{1}f,u\rangle _{M}=
\langle P_{s}\tilde{f},u\rangle =\langle\tilde{f},P_{-s-1}u\rangle
,$$
where $\tilde{f}$\ is an arbitrary continuation of $f$ on the
manifold $S.$ On the other hand
$$\left\langle r_{M}I_{1}^{\ast }I_{1}f,u\right\rangle _{M}=
\left\langle f,\left( r_{M}I_{1}^{\ast }I_{1}\right) ^{\ast
}u\right\rangle _{M}.$$ Since $\tilde{f}$ is arbitrary, then
equality $\langle \tilde{f},P_{-s-1}u\rangle =\langle f,\left(
r_{M}I_{1}^{\ast }I_{1}\right) ^{\ast }u\rangle _{M}$ implies
$\left( r_{M}I_{1}^{\ast }I_{1}\right) ^{\ast
}=r_{U_{1}}P_{-s-1}=r_{U_{1}}I_{1}^{\ast }I_{1}.$

Because of ellipticity the equality $r_{U_{1}}Pu =0$ implies
smoothness $u |_{U_{1}},$ and then $u \in H^{-s-1}\left( M\right)
$ implies of $u \in C_{0}^{\infty }( U_{1})$. Since $r_{U_{1}}Pu
=I_{1}^{\ast }I_{1}u ,$ then $I_{1}^{\ast }I_{1}u
=0\Longrightarrow I_{1}u =0\Longrightarrow u =0.$

\end{proof}

Now we are ready to prove Theorem 1.4.

\begin{proof}
Let  $I,I_{1}$ be the geodesic X-ray transforms on $M$ and $U_{1}$
respectively. From Theorem $3.1$ it follows that for any $h\in
C^{\infty }(M)$ there exists\  $f\in C^{\infty }(\overline{U})$,
such that $r_{M}I_{1}^{\ast }I_{1}f=h.$ Then $u^{f}\in C^{\infty
}(\Omega (\dot{N}))$ where $\dot{N}$ denotes the interior of $N$.
Let $w=2u_{+}^{f}|_{\partial _{+}\Omega \left( M \right)},$ where
$u_{+}^{f}$ is the even part with respect to $\xi$. Then it easy
to see, that $w_{\psi }=2u_{+}^{f}|_{\Omega(M)}$ and $I^{\ast
}w=h$. The function $w\in C_{\alpha }^{\infty }(\partial
_{+}\Omega (M))$ since $w_{\psi }\in C^{\infty }(\Omega (M))$.
\end{proof}

\bigskip

\section{\bf Scattering relation and folds}
\setcounter {equation}{0}

In this section we prove Lemma 1.1.

As indicated before we embed $(M,g)$ into a compact manifold
$(S,g)$ with no boundary. Let $(N,g)$ be an arbitrary neighborhood
in $(S,g)$ of the manifold $(M,g)$, such that any geodesic
$\gamma(x,\xi,t),(x,\xi)\in \Omega(N)$ intersects the boundary
$\partial N$ transversally. Then the length of the geodesic ray
$\tau$ is a smooth function on $\Omega(\dot{N})$ and the map
$\psi:\partial\Omega(M)\rightarrow\partial_{+}\Omega(N)$, defined
by
\begin{equation}
\phi(x,\xi)=\varphi_{\tau(x,\xi)}(x,\xi),\quad
(x,\xi)\in\partial\Omega(M),
\end{equation}
is smooth as well. Moreover it turns out $\phi$  is a fold map
with fold $\partial_{0}\Omega(M)$. This fact will be proved in the
next Theorem. Once this is proven Lemma 1.1 follows from [H]
Theorem C.4.4. From the assumption $A_{+}w\in
C^{\infty}(\partial\Omega(M))$ we deduce the existence of a
smooth function $v$ on a neighborhood of the range
$\phi(\partial\Omega(M))$ such that $w=v\circ\phi$. Consider
function $w_{\psi }=w\circ \alpha \circ \psi .$ Change notation
$\psi $\ to $\psi _{M},$ keeping $w_{\psi }.$ Denote by\ $\psi
_{N}$\ the map, analogical to $\psi _{M},$
$$\psi _{N}\left( x,\xi \right) =
\varphi _{\tau \left( x,\xi \right) }\left( x,\xi \right) ,
\;\;\left( x,\xi \right) \in \Omega \left( N\right) .$$ Then
$w_{\psi }=v\vspace{1pt}\circ \phi \circ \alpha \circ \psi _{M}.$
It easy to see, that $\phi \circ \alpha \circ \psi _{M}=\psi
_{N}|_{\Omega (M)}.$ Since\ the map\ $\psi _{N}$\ is smooth on
$\Omega \left( M\right) ,$ then $w_{\psi }\in C^{\infty }\left(
\Omega \left( M\right) \right) $, i.e. $w\in C_{\alpha }^{\infty
}( \partial _{+}\Omega (M)) .$ Thus Lemma 1.1 is proven once we
show that $\phi$ is a fold.

\begin{theorem}
Let $(M,g)$ be a strictly convex, non-trapping manifold and $N$ an
arbitrary neighborhood of $M$, such, that any geodesic
$\gamma(x,\xi, t),(x,\xi)\in \Omega(\dot{N})$ intersects the
boundary $\partial N$ transversally. Then the map $\phi$, defined
by (4.1) is a fold with fold $\partial_{0}\Omega(M)$.
\end{theorem}

First we recall the definition of a Whitney fold.
\begin{definition}
Let $M,N$ be  $C^\infty$ manifolds of the same dimension and let $f:M\longrightarrow N$
be a $C^\infty$ map with $f(m)=n.$ The function $f$ is a Whitney fold (with fold L) at $m$ if $f$ drops rank by one simply at $m$, so that $\{ x; df(x) \mbox
{ is singular } \}$ is
a smooth hypersurface near $m$ and $\mbox{ker}(df(m))$ is transverse to $T_{m}L.$
\end{definition}

We rewrite the definition of a fold below in a form that we will use.

Let $\gamma:(a,b)\rightarrow M, 0\in(a,b)$ be an arbitrary curve
in M. Then $f$ induces the curve $\gamma_{1}:(a,b)\rightarrow N$
in the manifold $N$, $\gamma_{1}=f\circ\gamma$. The derivative
$f'(m)$ of the map $f$ at the point $m$ maps the tangent vector
$\xi=\dot{\gamma}(0)\in T_{m}(M)$ into the tangent vector
$\eta=\dot{\gamma}_{1}(0)=f'(m)\xi\in T_{f(m}(N)$. The
(acceleration) vectors $\ddot{\gamma}(0), \ddot{\gamma}_{1}(0)$
belong to the spaces $T_{\xi}(T_{m}(M))$ and
$T_{\eta}(T_{f(m)}(N))$ respectively, which are tangent to the
vector spaces $T_{m}(M)$ and $T_{f(m)}(N)$ and may be identified
with them. Fix a point $m\in M$ and nonzero vector $\xi\in
Kerf'(m)$. Let us consider the map $L_{\xi}:T_{m}\rightarrow
T_{f(m)},\quad L_{\xi}(\ddot{\gamma}(0))=\ddot{\gamma}_{1}(0)$. A
calculation in local coordinates gives that the range of the map
$L_{\xi}$ is an affine subspaces in $T_{f(m)}(N),
RanL_{\xi}=a(\xi)+Ranf'(m)$, where $a(\xi)$ is a vector, depending
on $\xi$. Moreover the function
$$Hessf(m)(\xi,Y)=\langle L_{\xi}(X),Y\rangle, Y\in (Ranf'(m))^{\perp}$$
does not depend on $X$. Here $\langle,\rangle$ corresponds to the
duality of the spaces $T_{n}(N)$ and $T^{\ast}_{n}(N)$. It follows
immediately from the that the definition of a fold is equivalent
to
\begin{definition}
The map $f$ is a fold at the point $m$
if
\begin{enumerate}[(i)]
\item $\mbox{ dimKer }f'(m)=1.$
\item
$Hessf(m)(\xi,Y)\neq0.$
\end{enumerate}
\end{definition}

\begin{remark}
Condition 4.2 (i) implies that $\mbox{ dimCoker }f'(m)=1$ and
$\mbox{ Ran }L_{\xi}$ is an hyperplane. and condition 4.2(ii)
implies that the hyperplane $RanL_{\xi}$ does not pass trough the
origin.
\end{remark}

The following proposition is immediate.

\begin{proposition}
Let $f:M\rightarrow N$ be a smooth map, $\dim
M=\dim N$. Let $e:M_{0}\rightarrow M$ be an hypersurface and the map
$\ f\circ e:M_{0}\rightarrow f\left( M_{0}\right) $ is non-singular.
Then $\dim Ker\,f^{\prime }\left( x\right) \leq 1$ for any point
$x\in M_{0}$.
\end{proposition}

Now we start proving Theorem 4.1.

\begin{proof}

We first prove Condition 4.2(i). The smooth map
$\phi|_{\partial_{0}\Omega(M)}:\partial_{0}\Omega(M)\rightarrow
\phi(\partial_{0}\Omega(M))$ is bijective and therefore $\mbox{
dimKer }\phi'(x_{0},\xi_{0})\leq 1$ for any point
$(x_{0},\xi_{0})\in \partial_{0}\Omega(M)$. We prove that the
vector $X_{0}\in T_{(x_{0},\xi_{0})}(\partial\Omega(M)), X_{0 h
}=\xi_{0}, X_{0 v}=0$, where $X_{0 h},X_{0 v}$ are its horizontal
and vertical components (see section 2) belongs to $\mbox{ Ker
}\phi'(x_{0},\xi_{0})$. Represent $\phi=\psi_{N}\circ e$, where
$e:\partial\Omega(M)\rightarrow \Omega(N)$ is the embedding. Then
$\phi'(x,\xi)X=\psi_{N}'(x,\xi)X$ (we identify here the vector $X$
with its embedding $e'(x,\xi)X\in
T_{(x,\xi)}(\partial\Omega(M))$). Then $X_{0}\in Ker\,\phi
^{\prime }(x_{0},\xi _{0})$ iff $X_{0}\in T_{(x_{0},\xi _{0})}(
\partial \Omega ( M) )\bigcap Ker \psi
_{N} ^{\prime }(x_{0},\xi _{0}).$ The fact that  $X_{0}\in
T_{(x_{0},\xi _{0})}\left( \partial \Omega \left( M\right) \right)
$ follows from equality $\left( \nu \left( x_{0}\right) ,\xi
_{0}\right) =0$ (see  (2.2)). The equality \ $\psi _{N}^{\prime
}(x_{0},\xi _{0})X_{0}=(\nabla \psi _{N}(x_{0},\xi _{0}),\xi
_{0})=0$ is obvious, since for any point $(x,\xi)\in \Omega(N)$
the point $\psi_{N}(x,\xi)$ is left fixed by geodesic flow and
$(\nabla\psi _{N}(x,\xi),\xi)=\mathcal{H}\psi_{N}(x,\xi)=0$. Thus,
the vector $X_{0}=(\xi_{0},0)\in \mbox{ Ker }\phi'(x_{0},\xi_{0})$
and the first condition in the definition of a fold has been
verified.

Let us prove the second condition in the definition of a fold by
contradiction. Assume that
$\mbox{ Hess }\phi(x_{0},\xi_{0})(X_{0},Y^{0})=0$, where $Y^{0}\in
Ker(\phi'(x_{0},\xi_{0}))^\ast$. We will show then $Y^{0}=0$.
Denote $(y,\eta)=\phi(x,\xi)$. Then in local coordinates we have
for the covector $Y^{0}=(Y^{0}_{1},...,Y^{0}_{2n})$:
\begin{equation}
Y_{i}^{0}\nabla _{k}^{\parallel }y^{i}(x_{0},\xi_{0} )+
Y_{n+i}^{0}\nabla _{k}^{\parallel }\eta ^{i}(x_{0},\xi_{0} )=0,
\end{equation}
\begin{equation}
Y_{i}^{0}\partial_{k}y^{i}(x_{0},\xi_{0} )+
Y_{n+i}^{0}\partial_{k}\eta ^{i}(x_{0},\xi_{0})=0.
\end{equation}
We write down these equations using the horizontal and vertical
parts of the covector $Y^{0}$ (see section 2) and the Jacobi
fields $A_{(k)}\left( x_{0},\xi _{0},t\right) ,\;B_{\left(
k\right) }\left( x_{0},\xi _{0},t\right) ,\;k=1,...,n$\ on
the geodesic $\gamma (x_{0},\xi _{0},t):$
\[
 A_{(k)}=\nabla _{k}^{\parallel }\gamma +\dot{\gamma}\nabla
_{k}^{\parallel }\tau ,\quad B_{k}=\partial
_{k}\gamma+\dot{\gamma}\partial _{k}\tau.
\]
Then the equations (4.2),(4.3) can be written in the form
\[
\lbrack (A_{(k)},Y_{h}^{0})+(DA_{(k)},Y_{v}^{0})\rbrack |_{t=\tau
}=0,
\]
\[
\lbrack (B_{(k)},Y_{h}^{0})+(DB_{(k)},Y_{v}^{0})\rbrack |_{t=\tau
}=0,
\]
where $D$ means covariant derivative along geodesic. Without lost
of generality one can take $Y^{0}_{v}=Z^c|_{t=\tau },$\
$Y^{0}_{h}=-DZ^c|_{t=\tau }$, where $Z^c$ is some Jacobi field on
$\gamma(x_0,\xi_0, t)$. It is known that for any Jacobi fields
$X,Y$, the Wronskian $\{X,Y\}=(DX,Y)-(X,DY)$ is constant. Then we
get
\[
\{Z^{c},A_{(k)}\}|_{t=0}=0,\quad \{Z^{c},B_{(k)}\}|_{t=0}=0.
\]
The initial data is given by
\[
A_{(k)}^{j}|_{t=0}=\delta _{k}^{j}-\nu _{0}^{i}\nu _{0k}+\xi
_{0}^{j}\nabla _{k}^{\parallel }\tau ,\quad DA_{(k)}|_{t=0}=0,
\]
\[
B_{(k)}|_{t=0}=\xi _{0}\partial _{k}\tau ,\quad
(DB_{(k)})^{j}|_{t=0}=\delta _{k}^{j}-\xi _{0}^{j}\xi _{0k},
\]
where $\nu_{0}=\nu(x_{0}).$ We have
\[
Z^c|_{t=0}=a\xi_{0}+c\partial\tau,\quad
DZ^c|_{t=0}=b\nu_{0}-c\nabla^{\parallel}\tau,
\]
where $a,b,c$ are some constants, depending on the point
$(x_{0},\xi_{0})$. From these general solutions we find one
satisfying the  condition that $Y^{0}\in
T_{\phi(x,\xi)}(\partial\Omega(N))$, or the conditions (2.2),
\begin{equation}
(Y^{0}_{h},\nu_{N}(y))=-(DZ^c|_{t=\tau},\nu_{N}(y))=0,
\end{equation}
\begin{equation}
(Y^{0}_{v},\eta)=(Z^c|_{t=\tau},\eta)=0,
\end{equation}
where $\nu_{N}$ is the normal to $\partial N$. Consider at first
the condition (4.5). Using the Jacobi field
$J=(\tau-t)\dot{\gamma}$ it can be rewritten in the form
\[
\{Z^c,J\}|_{t=\tau}=\{Z^c,J\}|_{t=0}=0,
\]
and, therefore, $a=-c\tau$. Thus,
\[
Z^c|_{t=0}=c(\partial\tau-\xi_{0}\tau),\quad
DZ^c|_{t=0}=b\nu_{0}-c\nabla^{\parallel}\tau.
\]
Now consider our hypothesis that $\mbox{ Hess
}\phi_{(x_{0},\xi_{0})}(X_{0},Y^{0})=0$. A straightforward
calculation shows that in local coordinates we have
\[
\mbox{ Hess
}\phi_{(x_{0},\xi_{0})}(X_{0},Y^{0})=-B(x_{0},\xi_{0})\{Z,P\}|_{t=\tau}=
-B(x_{0},\xi_{0})\{Z^{c},P\}|_{t=0},
\]
where $B=-\xi^{i}\xi^{j}\nabla^{\parallel}_{i}\nu_{j}$ is the
second fundamental form of the boundary $\partial M$ and
$$P=(\nu,\nabla\gamma)+(\nu ,\nabla \tau )\dot{\gamma}$$ is a
Jacobi field. Since $B>0$, our hypothesis gives that
$\{Z,P\}|_{t=0}=0$. Because of the initial data
$P(x_{0},\xi_{0},0)=\nu_{0}+(\nu_{0} ,\nabla \tau )\xi, \quad
DP(x_{0},\xi_{0},0)=0$ we obtain $b=-c(\nu_{0},\nabla\tau)$ and
\[
Z^c|_{t=0}=c(\partial\tau-\xi_{0}\tau),\quad
DZ^c|_{t=0}=-c\nabla\tau.
\]

Now we prove that for any point $( x,\xi) \in \Omega ( \dot{N}) $
the vector field
\begin{equation}
Z\left( x,\xi ,t\right) =(\partial \tau -\xi \tau )\left( \varphi
_{t}\left( x,\xi \right) \right) ,\;t\in \lbrack -\tau \left(
x,-\xi \right) ,\tau \left( x,\xi \right) \rbrack
\end{equation}
on the geodesic $\gamma \left( x,\xi ,t\right) $ is a Jacobi field and
\begin{equation}
DZ\left( x,\xi ,t\right) =\left( \nabla \tau \right) \left(
\varphi _{t}\left( x,\xi \right) \right).
\end{equation}
Applying to the equation $\mathcal{H}\tau =-1$ the operator
$\mathcal{H}\partial $ and using (2.3)-(2.5) we get
$\mathcal{H}^{2}X+RX=0,\;\;\left( RX\right) ^{i}(x,\xi
)=R_{jkl}^{i}\left( x\right) \xi ^{j}X^{k}\xi ^{l},$ where the
semibasic vector field $X=\partial \tau -\xi \tau .$ The operators
$\mathcal H$ and $D$ are related by:
\begin{equation}
\left( \mathcal{H} X\right)\circ \varphi _{t} =D(X\circ \varphi
_{t}) ,
\end{equation}
where $X$ is an arbitrary semibasic vector field. Then (4.6) is a
Jacobi field. Since from $\mathcal{H}\tau =-1$ it follows
$\mathcal{H}(
\partial \tau -\xi \tau ) =-\nabla \tau ,$ we obtain (4.7) from
(4.8).

Thus, $Z^{c}( x_{0},\xi _{0},t) =cZ( x_{0},\xi _{0},t) $ and
condition (4.4) gives
\begin{equation}
0=(DZ|_{t=\tau },\nu_{N} (y))=-c(\nabla \tau
( y,\eta ) ,\nu _{N}(y)).
\end{equation}

Now we show that
\begin{equation}
\nabla \tau |_{\partial _{-}\Omega \left( N\right) }=-\frac{\nu
_{N}}{\left( \nu _{N},\xi \right) }.
\end{equation}
Let $h$\ be a\ smooth\ function on $N$, such that $h|_{\partial
M}=0,\;\nabla h|_{M}\neq 0.$ Then since $h(x,\xi ,\tau (x,\xi
))=0,\;(x,\xi )\in \Omega (N)$, we have that
\[
\nu _{i}^{N}(y)\nabla \gamma _{(k)}^{i}(x,\xi ,\tau (x,\xi ))+
\left( \nu _{N}\left( y\right) ,\eta \right) \nabla _{k}\tau
(x,\xi )=0,
\]
where $\left( y,\eta \right) =\varphi _{\tau \left( x,\xi \right)
}(x,\xi )\in \partial _{-}\Omega (N)$. Change here $(x,\xi )$ to
$\varphi _{t}(x,\xi ),\;t\in \left[ 0,\tau \left( x,\xi \right)
\right] .$ Then the point $\left( y,\eta \right) $ does not
change. Using the initial data $\nabla \gamma _{(k)}^{i}(x,\xi
,0)=\delta _{k}^{i}$ we get in the limit $\nu _{N}(y)+(\nu
_{N}(y),\eta )\nabla \tau (y,\eta )=0.$ Thus, (4.10) is proved and
we conclude from (4.9) that $c=0.$

Thus, $\mbox{ Hess }\phi_{(x_{0},\xi_{0})}(X_{0},Y^{0})=0$ iff
$\mbox{ Ker }(\phi'(x_{0},\xi_{0}))^\ast=\{0\}$. We have
finished the proof of the Theorem.

\end{proof}

\bigskip

\section{\bf The Hilbert transform and geodesic flow}
\setcounter {equation}{0}
\bigskip

In this section we prove Theorem 1.5 in the introduction.


\

Let $H$ be the Hilbert transform as defined in (1.3)

We have that
$H$ is a unitary operator in the space
$L^{2}_{0}(\Omega_{x})=\{u\in L^{2}(\Omega_{x}): u_{0}=0\}$,
$$(u,v)=(Hu,Hv),\quad \forall u,v\in L^{2}_{0}(\Omega_{x}),$$
$$H^{2}(u)=-u,\quad \forall u\in L^{2}_{0}(\Omega_{x}).$$
Clearly, all these properties remain the same if we change $\Omega_{x}$
to $\Omega(M)$. \



In order to prove Theorem 1.4 we need the following commutator formula
which is valid for Riemannian manifolds of any dimension
\begin{lemma}
Let $u$ be a smooth function on the manifold $\Omega
^{2}(M)=\bigcup_{x\in M}\Omega^{2}_{x},\quad
\Omega^{2}_{x}=\{(x,\xi,\eta):\xi,\eta\in \Omega_{x}\}.$ Then
\begin{equation}
\nabla \int\limits_{\Omega _{x}}u(x,\xi ,\eta)
d\Omega _{x}\left( \eta \right) = \int\limits_{\Omega _{x}}\nabla^{(2)}
u(x,\xi ,\eta ) d\Omega _{x}\left( \eta \right),
\end{equation}
where $\nabla^{(2)}$ under the integral sign in (5.1) denotes the
horizontal derivative on $\Omega ^{2}(M)$,
\end{lemma}
\[
\nabla^{(2)}_{j} u(x,\xi ,\eta )=(\frac{\partial}{\partial
x^{j}}-\Gamma^{i}_{jk}\xi^{k}\partial_{i(\xi)}-\Gamma^{i}_{jk}\eta^{k}\partial_{i(\eta)})u(x,\xi
,\eta ).
\]
Notice that the horizontal tangential derivative can be defined on
$T(M)\times T(M)$ in a similar fashion to the case of $T(M)$ in
section 2.

\begin{proof}
Let $\varphi \in C_{0}^{\infty }\left( \R^{+}\right) $
be arbitrary function. We define the function $v$ on $T^2\left( M\right)
$ by
\[
v(x,\xi ,\eta )=\varphi \left( \left| \eta \right|
\right) u(x,\xi /\left| \xi \right| ,\eta /\left| \eta \right| )
\]
Let us consider the integral
\[
S(x,\xi )=\int\limits_{T_{x}(M)}v(x,\xi ,\eta )dT_{x}\left( \eta \right) .
\]
Identifying $T_{x}(M)$ with $R^{n}$ we have
\[
S(x,\xi )=\int\limits_{R^{n}}v(x,\xi ,\eta )\sqrt{\det g\left( x\right) }d\eta .
\]
Then
\[
\nabla _{j}S=\frac{\partial S}{\partial x^{j}}- \Gamma
_{jk}^{i}\xi ^{k}\frac{\partial S}{\partial \xi ^{i}}=
\int\limits_{R^{n}}(\frac{\partial v}{\partial x^{j}}- \Gamma
_{jk}^{i}\xi ^{k}\frac{\partial v}{\partial \xi ^{i}})\sqrt{\det
g}d\eta +\int\limits_{R^{n}}v\frac{\partial \ln \sqrt{\det g\left(
x\right) }}{\partial x^{j}}\sqrt{\det g}d\eta.
\]
Since $\partial \ln \sqrt{\det g}/dx^{j}=\Gamma _{jk}^{k}$ we
rewrite the last integral in the form
$$\int\limits_{R^{n}}v\frac{\partial }{\partial \eta ^{k}}
\left( \Gamma _{jl}^{k}\eta ^{l}\right) \sqrt{\det g}d\eta .$$
Then
$$\nabla _{j}S=\int\limits_{R^{n}}(\frac{\partial v}{\partial x^{j}}-
\Gamma _{jk}^{i}\xi ^{k}\frac{\partial v}{\partial \xi ^{i}}-
\Gamma _{jl}^{k}\eta ^{l}\frac{\partial v}{\partial \eta
^{k}})\sqrt{\det g}d\eta.$$ Since $$(\frac{\partial }{\partial
x^{j}}- \Gamma _{jk}^{i}\xi ^{k}\frac{\partial }{\partial \xi
^{i}}- \Gamma _{jl}^{k}\eta ^{l}\frac{\partial }{\partial \eta
^{k}})\left| \eta \right| =0,$$ then after changing to spherical
coordinates we obtain
\begin{equation}
\nabla S(x,\xi )=\int\limits_{0}^{\infty }
\varphi \left( t\right) t^{n-1}dt \int\limits_{\Omega _{x}}\nabla
u(x,\xi ,\eta ) d\Omega _{x}\left( \eta \right) .
\end{equation}

Now $S$ in spherical coordinates is given by
\begin{equation}
S(x,\xi )=\int\limits_{0}^{\infty }\varphi \left( t\right) t^{n-1}dt
\int\limits_{\Omega _{x}}u(x,\xi ,\eta ) d\Omega _{x}\left( \eta
\right)
\end{equation}

We conclude (5.1) using (5.2),(5.3).

\end{proof}

\bigskip

Now we prove Theorem 1.4.

\begin{proof}

A straightforward calculation gives that
$$\nabla
\frac{1+\left( \xi ,\eta \right) }{\left( \xi _{\perp },\eta
\right) }=0$$

and therefore we have
\[
\nabla H u(x,\xi )=\frac{1}{2\pi }\int\limits_{\Omega
_{x}}\frac{1+\left( \xi ,\eta \right) }{\left( \xi _{\perp },\eta
\right) }\nabla u(x,\eta )d\Omega _{x}\left( \eta \right).
\]

For any pair vectors $\xi ,\,\eta \in \Omega _{x}$ we have
$$\eta=(\xi ,\eta )\xi +(\xi _{\perp },\eta )\xi _{\perp },  \eta
_{\perp }=-(\xi _{\perp },\eta )\xi +(\xi ,\eta )\xi _{\perp
},\quad \left( \xi ,\eta \right) ^{2}+\left( \xi _{\perp },\eta
\right) ^{2}= 1.$$

Then $$\eta \frac{1+\left( \xi ,\eta \right) }{\left( \xi _{\perp
},\eta \right) }=\xi \frac{\left( \xi ,\eta \right) +\left( \xi
,\eta \right) ^{2}}{\left( \xi _{\perp },\eta \right) }+\xi
_{\perp }(1+(\xi ,\eta ))=$$ $$=\xi \frac{\left( \xi ,\eta \right)
+1}{\left( \xi _{\perp },\eta \right) }-\xi \left( \xi _{\perp
},\eta \right) +\xi _{\perp }(\xi ,\eta )+\xi _{\perp }=\xi
\frac{1+\left( \xi ,\eta \right) }{\left( \xi _{\perp },\eta
\right) }+\xi _{\perp }+\eta _{\perp }.$$

Thus \[
H\mathcal{H}u=\mathcal{H}Hu+\mathcal{H}_{\perp
}u_{0}+(\mathcal{H}_{\perp }u)_{0}
\]
and
Theorem 1.4 is proved.

\end{proof}

\bigskip



\end{document}